\documentclass{amsart}
\usepackage{amssymb,amsmath}

\newtheorem{theorem}{Theorem}
\newtheorem{lemma}{Lemma}

\newcommand{\bt}{\begin{theorem}}
\newcommand{\et}{\end{theorem}}
\newcommand{\bl}{\begin{lemma}}
\newcommand{\el}{\end{lemma}}
\newcommand{\beal}{\begin{align*}}
\newcommand{\enal}{\end{align*}}
\newcommand{\bq}{\begin{eqnarray*}}
\newcommand{\eq}{\end{eqnarray*}}
\newcommand{\be}{\begin{eqnarray}}
\newcommand{\ee}{\end{eqnarray}}
\newcommand{\beq}{\begin{equation}}
\newcommand{\eeq}{\end{equation}}
\newcommand{\benum}{\begin{enumerate}}
\newcommand{\eenum}{\end{enumerate}}
\newcommand{\ba}{\begin{array}}
\newcommand{\ea}{\end{array}}

\begin{document}
\title{On distinct consecutive differences}
\author{J\'ozsef Solymosi}
\thanks{This research was supported by NSERC and OTKA grants.}
\address{Department of Mathematics\\
University of British Columbia\\
Vancouver} \email{solymosi@math.ubc.ca}

\begin{abstract}
We show that if $A=\{a_1,a_2,\ldots ,a_k\}$ is a monotone
increasing set of numbers, and the differences of the consecutive
elements are all distinct, then $|A+B|\geq c|A|^{1/2}|B|$ for any
finite set of numbers $B$. The bound is tight up to the constant
multiplier.
\end{abstract}

\maketitle

\section{introduction}
\noindent Given two sets of numbers, $A$ and $B$, the {\em sumset}
of $A$ and $B$, denoted by $A+B,$ is
\[
A+B = \{a+b: a,b \in A \text{ and } b \in B\}
\]
Let $A=\{a_1,a_2,\ldots ,a_k\}$ be a finite set of real numbers
with the property that
\begin{equation}
a_i-a_{i-1}<a_{i+1}-a_i
\end{equation}
for any $1<i<k.$ Sets with this property are said to be {\em
convex} sets. Answering a question of Erd\H os, Hegyv\'ari
\cite{HE} proved that if $A$ is convex then $|A+A|\geq ck\log
k/\log\log k.$

Hegyv\'ari's result was later improved by Elekes, Nathanson, and
Ruzsa \cite{ENR}. They proved that if $A$ is convex, then
$|A+B|\geq ck^{3/2}$ for any set $B$ with $|B|=k.$ In this paper
we extend this result for sets with distinct consecutive
differences. Set $A$ has distinct consecutive differences if for
any $1\leq i,j\leq k,$ $a_{i+1}-a_i=a_{j+1}-a_j$ implies $i=j.$

\bt Let $A$ and $B$ be finite sets of real numbers with $|A| = k$
and $|B| = \ell.$ If $A$ has distinct consecutive differences,
then
\[
|A+B| \geq \frac{k\sqrt{\ell}}{3}.
\]
In particular, if $k = \ell,$ then
\[
|A+B| \geq \frac{k^{3/2}}{3}  .
\]
\et

The basic idea behind the proof is the following. The sumset $A+B$
consists of $|B|$ translates of $A$. The translates of two
consecutive elements of $A$ are typically not "far" from each
other in the sumset $A+B.$ Also, from a translate of two
consecutive elements, $b+a_i$,$b+a_{i+1}$ we can recover the value
of $b$, since all of the consecutive differences are distinct.
Then the number of "close" pairs in $A+B$ should be large, around
$|A||B|$, therefore $A+B$ is also large.

\medskip

In the second part of the paper we extend the result for two sets.
As an application we show that for any convex function $F$, and
finite sets of real numbers, $A,B,$ and $C$, if $|A|=|B|=|C|=n,$
then $\max\{|A+B|,|F(A)+C|\geq cn^{5/4}.$

\medskip

Along the lines of the proof one can prove the "statistical"
version of Theorem 1 and Theorem 3. We state the analogue of
Theorem 1 without working out the details of the proof.

\bt Let $A=\{a_1,a_2,\ldots ,a_n\}$ is a monotone increasing set
of numbers. If the difference set of the consecutive elements,
$D=\{a_{i+1}-a_i : 1\leq i \leq n-1\}$, is large, $|D|\geq \delta
|A|$, then $|A+B|\geq c|A|^{1/2}|B|$ for any finite set of numbers
$B$, where $c$ depends on $\delta$ only. \et

\section{Distinct consecutive differences}

\proof {\em of Theorem 1.} Since $|A+B| \geq \min(k,\ell),$ the
result is immediate for $\min(k,\ell) \leq 3,$ so we can assume
that $k \geq 3$ and $\ell \geq 3.$

Let $A = \{a_1,a_2,\ldots,a_k\}$ and $B = \{ b_1,b_2,\ldots, b_{\ell}\}$,
where
$a_1 < a_2 < \cdots < a_k$ and $b_1 < b_2 < \cdots < b_{\ell}$.
Let
\[
A+B = C = \{c_1,c_2,\ldots,c_m\},
\]
where $c_1< c_2 < \cdots < c_m$ and $m = |A+B|.$

Let $1 \leq i \leq k-1$ and $1 \leq j \leq \ell.$ A {\em pair} is
a two-element subset of $C$ of the form \beq \label{soly:pair}
\{a_i+b_j, a_{i+1}+b_j\}. \eeq Suppose that $\{c,c'\}$ is a pair
and $c < c'$. Since the set $A$ has distinct consecutive
differences, there is a unique integer $i$ such that
\[
c'-c = a_{i+1}-a_i.
\]
It follows that there is a unique integer $j$
such that
\[
c - a_i = c' - a_{i+1} = b_j.
\]
Therefore, if
\[
\{a_i+b_j, a_{i+1}+b_j\} = \{a_{i'}+b_{j'}, a_{i+1}+b_{j'}\},
\]
then $i = i'$ and $j = j'$, and so the sumset $C$ contains
exactly $(k-1)\ell$ pairs.

Let $s_0,s_1,\ldots, s_t$ be integers such that
\[
0 = s_0 < s_1 < s_2 < \cdots < s_{t-1} < s_t = m.
\]
We partition the set $C$ into $t$ pairwise disjoint sets $C_1,\ldots, C_t$ as follows:
\[
C_1 = \{c_{1},\ldots,c_{s_1}\},
\]
\[
C_2 = \{c_{s_1 +1},\ldots,c_{s_2}\},
\]
and, in general, for $u = 1,\ldots, t,$
\[
C_u = \{c_{s_{u-1}+1},\ldots,c_{s_u}\}.
\]
Then
\[
|C_u| = s_u - s_{u-1} \quad\text{ for $u = 1,\ldots, t$}.
\]
Fix an integer $j \in \{1,2,\ldots,\ell\}$, and consider the increasing sequence
\[
a_1+b_j < a_2+b_j < \cdots < a_k+b_j.
\]
Let $k_{j,u}$ denote the number of elements of this sequence that belong to the set $C_u.$
Then
\[
\sum_{u=1}^t k_{j,u} = k.
\]
If $k_{j,u} \geq 1,$ then the set $C_u$ contains exactly $k_{j,u}-1$ pairs
of the form~(\ref{soly:pair}),
and so the number of pairs with fixed $j$
contained in the sets $C_1,\ldots,C_t$ is
\[
\sum_{u=1\atop k_{j,u} \geq 1}^t (k_{j,u}-1) = k - \sum_{u=1\atop k_{j,u} \geq 1}^t 1\geq k-t.
\]
Since the set $C$ contains exactly $(k-1)\ell$ distinct pairs,
it follows that the total number of pairs contained in the sets $C_1,\ldots,C_t$ is at least
\[
\sum_{j=1}^{\ell}  \sum_{u=1\atop k_{j,u} \geq 1}^t (k_{j,u}-1)
\geq \ell(k-t).
\]
We can obtain a simple upper bound for the total number of pairs contained in the sets $C_1,\ldots,C_t$ as follows:
For $u = 1,\ldots,t,$ the set $C_u$ contains at most $|C_u| \choose 2$ pairs,
and so the number of pairs contained in the sets $C_1,\ldots,C_t$ is at most
\[
\sum_{u=1}^t {|C_u| \choose 2}.
\]
Therefore,
\[
\sum_{u=1}^t {|C_u| \choose 2} \geq \ell(k-t).
\]

We specialize this inequality as follows:
Let
\[
t = \left[ \frac{k}{2}\right]
\]
and
\[
m = qt+r,
\]
where
\[
0 \leq r \leq t-1.
\]
Then
\[
q \leq \frac{m}{t} \leq  \frac{2m}{k-1}.
\]
Choose the integers $s_1,\ldots,s_{t-1}$ such that
\[
|C_u| = q+1 \qquad\text{for $u=1,\ldots,r$}
\]
and
\[
|C_u| = q\qquad\text{for $u=r+1,\ldots,t$}.
\]
Then
\[
\sum_{u=1}^t {|C_u| \choose 2} \leq t{q+1\choose 2}
\leq \frac{k}{2}{\frac{2m}{k-1}+1\choose 2}
< \frac{k}{4}\left( \frac{2m}{k-1}+1\right)^2
\]
and so
\[
\frac{k}{4}\left( \frac{2m}{k-1}+1\right)^2
>  \frac{\ell k}{2}.
\]
This implies that
\[
m > \frac{k-1}{2}(\sqrt{2\ell} -1).
\]
If $k \geq 3$ and $\ell \geq 3,$ then
\[
m > \frac{k\sqrt{\ell}}{3}.
\]
This completes the proof.

\section{Distinct pairs of consecutive differences}

Let $A = \{a_1, a_2,\ldots, a_k\}$ and $A' = \{a'_1, a'_2, \ldots,
a'_k\}$ be nonempty sets of real numbers, where $a_1 < a_2 <
\cdots < a_k$ and $a'_1 < a'_2 < \cdots < a'_{k}$. Let $d_i =
a_{i+1}-a_i$ for $i = 1,\ldots,k-1$ and  $d'_i = a'_{i+1}-a'_i$
for $i = 1,\ldots,k-1.$ The sets $A$ and $A'$ have {\em distinct
pairs of consecutive differences} if there exists a one-to-one map
$\sigma: \{1,2,\ldots, k-1\} \rightarrow \{1,2,\ldots, k-1\}$ such
that the $k-1$ ordered pairs $\left(d_i,d'_{\sigma(i)}\right)$ are
distinct.

\bt Let $A$ and $A'$ be nonempty finite sets of real numbers such
that $k = |A| = |A'|$ and the sets $A$ and $A'$ have distinct
pairs of consecutive differences. Let $B$, and $B'$ be nonempty
finite sets of real numbers with $|B|=\ell$, and $|B'| = \ell'$
Then
\[
|A+B|\cdot|A'+B'| \gg \left( k^3\ell\ell'\right)^{1/2}.
\]
If $\ell = \ell' = k,$ then
\[
|A+B|\cdot|A'+B'| \gg k^{5/2}.
\]
\et

\proof Let $A = \{a_1, a_2,\ldots, a_k\}$ and $A' = \{a'_1, a'_2,
\ldots, a'_k\}$, where $a_1 < a_2 < \cdots < a_k$ and $a'_1 < a'_2
< \cdots < a'_{k}$. Let $d_i = a_{i+1}-a_i$ for $i = 1,\ldots,k-1$
and $d'_i = a'_{i+1}-a'_i$ for $i = 1,\ldots,k-1$. Let $\sigma:
\{1,2,\ldots, k-1\} \rightarrow \{1,2,\ldots, k-1\}$ be a
one-to-one map such that the $k-1$ ordered pairs
$(d_i,d'_{\sigma(i)})$ are distinct. Let $B = \{b_1, b_2,\ldots,
b_{\ell}\}$ and $B' = \{b_1, b_2, \ldots, b_{\ell'}\}$, where $b_1
< b_2 < \cdots < b_{\ell}$ and $b'_1 < b'_2 < \cdots <
b'_{\ell'}$.

Let $1 \leq i \leq k-1$, $1 \leq j \leq \ell,$ and $1 \leq j' \leq \ell'.$
We consider quadruples of the form
\beq  \label{soly:quad}
(a_i+b_j,a_{i+1}+b_j, a'_{\sigma(i)}+b'_{j'}, a'_{\sigma(i)+1}+b'_{j'}).
\eeq
Suppose that $1 \leq u \leq k-1$, $1 \leq v \leq \ell,$ and $1 \leq v' \leq \ell',$
and that
\[
(a_i+b_j,a_{i+1}+b_j, a'_{\sigma(i)}+b'_{j'}, a'_{\sigma(i)+1}+b'_{j'})
= (a_u+b_v, a_{u+1}+b_v, a'_{\sigma(u)}+b'_{v'}, a'_{\sigma(u)+1}+b'_{v'}).
\]
Then
\[
d_i = (a_{i+1}+b_j) - (a_i+b_j) = (a_{u+1}+b_v) - (a_u+b_v) = d_u
\]
and
\[
d'_{\sigma(i)} = (a'_{\sigma(i)+1}+b'_{j'}) - (a'_{\sigma(i)}+b'_{j'})
= (a'_{\sigma(u)+1}+b'_{v'}) - (a'_{\sigma(u)}+b'_{v'}) = d'_{\sigma(u)}.
\]
Therefore,
\[
(d_i,d'_{\sigma(i)}) = (d_u,d'_{\sigma(u)}).
\]
The sets $A$ and $A'$ have matching consecutive differences
with permutation $\sigma$, and so $i = u.$
Since $a_i+b_j = a_u+b_v = a_i + b_v,$ it follows that $b_j = b_v$ and $j = v$.
Similarly, $j' = v'.$
This implies that there are $(k-1)\ell\ell'$ distinct quadruples of the form~(\ref{soly:quad}).

Consider the sumsets $A+B$ and $A'+B'$.
\[
A+B = C = \{c_1,\ldots,c_m\},
\]
where $|C| = m$ and $c_1 < \cdots < c_m.$
Let
\[
A'+B' = C' = \{c'_1,\ldots,c'_{m'}\},
\]
where $|C'| = m'$ and $c'_1 < \cdots < c'_{m'}.$

Let $s_0,s_1,\ldots, s_t$ be integers such that
\[
0 = s_0 < s_1 < s_2 < \cdots < s_{t-1} < s_t = m.
\]
We partition the set $C$ into $t$ pairwise disjoint sets $C_1,\ldots, C_t$ as follows:
\[
C_1 = \{c_{1},\ldots,c_{s_1}\},
\]
\[
C_2 = \{c_{s_1 +1},\ldots,c_{s_2}\},
\]
and, in general, for $u = 1,\ldots, t,$
\[
C_u = \{c_{s_{u-1}+1},\ldots,c_{s_u}\}.
\]
Then
\[
|C_u| = s_u - s_{u-1} \quad\text{ for $u = 1,\ldots, t$}.
\]
Similarly, let $s'_0,s'_1,\ldots, s'_{t'}$ be integers such that
\[
0 = s'_0 < s'_1 < s'_2 < \cdots < s'_{t-1} < s'_t = m',
\]
and partition the set $C'$ into $t'$ pairwise disjoint sets $C'_1,\ldots, C'_{t'}$ as follows:
\[
C'_{u'} = \{ c'_{s_{u'-1}+1},\ldots,c'_{s_{u'}} \}
\]
for $u' = 1,\ldots, t'.$

Fix an integer $j \in \{1,2,\ldots,\ell\}$, and consider the increasing sequence
\[
a_1+b_j < a_2+b_j < \cdots < a_k+b_j.
\]
These $k$ numbers belong to sumset $C.$
Let $k_{j,u}$ denote the number of elements of this sequence that belong
to the set $C_u.$  Then
\[
\sum_{u=1}^t k_{j,u} = k.
\]
If $k_{j,u} \geq 1,$ then the set $C_u$ contains exactly $k_{j,u}-1$ subsets
of the form
\beq \label{soly:pair2}
\{a_i+b_j, a_{i+1}+b_j\},
\eeq
and so the number of pairs with fixed $j$ contained in the sets $C_1,\ldots,C_t$ is
\[
\sum_{u=1\atop k_{j,u} \geq 1}^t (k_{j,u}-1)
= k - \sum_{u=1\atop k_{j,u} \geq 1}^t 1\geq k-t.
\]
There are $k-1$ pairs of the form~(\ref{soly:pair2}), and so the number of such pairs that
do not belong to one of the sets $C_1,\ldots,C_t$ is at most $t-1$.
Therefore, for each $j = 1,\ldots,\ell,$ the number of quadruples of the form~(\ref{soly:quad}) whose first two coordinates do not belong to one of the sets $C_1,\ldots,C_t$ is at most
\[
(t-1)\ell'.
\]
Summing over $j$, we conclude that the number of quadruples of the form~(\ref{soly:quad}) whose first two coordinates do not belong to one of the sets $C_1,\ldots,C_t$ is at most
\[
(t-1)\ell\ell'.
\]
Similarly, the number of quadruples of the form~(\ref{soly:quad}) whose third and fourth coordinates do not belong to one of the sets $C'_1,\ldots,C'_{t'}$ is at most
\[
(t'-1)\ell\ell'.
\]
Therefore, the number of quadruples of the form~(\ref{soly:quad}) whose first pair
of coordinates does not belong to a set $C_u$ or whose last pair of coordinates does not belong to a set $C'_{u'}$ is at most $(t+t'-2)\ell\ell'.$
It follows that the number of quadruples of the form~(\ref{soly:quad}) whose first pair
of coordinates belongs to one of sets $C_u$ and whose last pair of coordinates also belongs to one of the sets $C'_{u'}$ is bounded below by
\[
(k-1)\ell\ell' - (t+t'-2)\ell\ell' = (k-t-t'+1)\ell\ell'.
\]

On the other hand, for each $u = 1,\ldots, t$ and $u' = 1,\ldots, t',$ the number of quadruples $(x,y,z,w)$ such that $x\neq y$ and $\{x,y\} \subseteq C_u$ for $u = 1,\ldots, t$,
and also $z\neq w$ and $\{z,w\} \subseteq C'_{u'}$ is exactly
\[
{|C_u| \choose 2}{|C'_{u'}| \choose 2}
\]
It follows that a simple upper bound for the number of quadruples of the form~(\ref{soly:quad}) whose first pair of coordinates belongs to one of sets $C_u$ and whose last pair of coordinates also belongs to one of the sets $C'_{u'}$ is
\[
\sum_{u=1}^t\sum_{u'=1}^{t'} {|C_u| \choose 2}{|C'_{u'}| \choose 2}
\leq \frac{1}{4}\sum_{u=1}^t\sum_{u'=1}^{t'} |C_u|^2|C'_{u'}|^2.
\]
Therefore,
\beq   \label{soly:est}
(k-t-t'+1)\ell\ell' \leq \frac{1}{4}\sum_{u=1}^t\sum_{u'=1}^{t'} |C_u|^2|C'_{u'}|^2.
\eeq

If we let $t=t'=k/4$ and $|C_u| = m/t = 4m/k$ and $C'_{u'}| =
4m'/k,$ then
\[
\frac{k\ell\ell'}{2} \leq \frac{1}{4}\frac{k^2}{16}\frac{16m^2}{k^2}\frac{16m'^2}{k^2} =
\frac{4(mm')^2}{k^2},
\]
and so
\[
(mm')^2 \geq \frac{k^3\ell\ell'}{8}.
\]
If $\ell = \ell' = k,$ then
\[
|A+B|\cdot|A'+B'| = mm' \gg k^{5/2}.
\]

\section{Remarks}

A simple consequence of Theorem 3 is the following result, which
was first proved by Elekes, Nathanson, and Ruzsa \cite{ENR}.

\bt For any strictly convex real function $F$, and finite sets of
real numbers, $A,B,$ and $C$, if $|A|=|B|=|C|=n,$ then
$\max\{|A+B|,|F(A)+C|\geq cn^{5/4}.$ In particular, $|A+F(A)|\geq
cn^{5/4}.$

\et

\proof

The two sets $A$ and $F(A)$ have distinct pairs of consecutive
differences with $\sigma = I,$ the identity map. For the second
inequality set $B=F(A)$ and $C=A.$

\medskip

A construction of Ruzsa \cite{R} shows that the bounds in Theorem
1 and 3 are tight. However, as we will see, in his construction
$A$ and $B$ have very different structures. It is possible, that
if $A=B$ in Theorem 1, then $|A+B|$ is much larger, maybe close to
$|A|^2.$ Here we sketch Ruzsa's construction. Let $S$ be set such
that all the differences in $S$ are distinct, $|S-S|={|S|\choose
2},$ and $|S|$ is odd. Then there is a list $L$ of the elements of
$S$ with repetitions, consists of $k={|S|\choose 2}$ elements,
such that the consecutive elements have distinct differences.
($L=(s_1,s_2,\ldots , s_k)$ where $s_{i+1}-s_i=s_{j+1}-s_j$
implies that $i=j.$) Now we are ready to define $A.$

$$A=\{i+s_i : 1\leq i \leq k\}$$

The set $A$ has the property that the consecutive differences are
distinct, and it is not difficult to see that $|A+[k]|\leq
|S|^3\leq  ck^{3/2},$ where $[k]$ denotes the first $k$ natural
numbers.

The example shows that Theorem 1 is sharp, however Erd\H os'
original question is still wide open; What is the smallest
possible size of $A+A$ if $A$ is a convex set of integers?

\medskip

{\bf Acknowledgement:} I thank Mel Nathanson for his help on
writing up the paper.


\begin{thebibliography}{10}

\bibitem{E}
Gy.~Elekes, On the number of sums and products. {\em Acta Arith.}
81 (1997), no. 4, 365--367.

\bibitem{EL}
Gy.~Elekes, Sums versus products in number theory, algebra and
Erd\H os geometry. in: Paul Erd\H os and his Mathematics. II,
Budapest, Bolyai Society Mathematical Studies, 11. (2002), .

\bibitem{ENR}
Elekes, M.~Nathanson, and I,~Ruzsa, Convexity and sumsets. {\em J.
Number Theory} 83 (2000), no. 2, 194--201.

\bibitem{R}
I.~Ruzsa, personal communication (2003)

\bibitem{HE}
N.~Hegyv\'ari, On consecutive sums in sequences. {\em Acta Math.
Hungar.} 48 (1986), no. 1-2, 193--200.




\end{thebibliography}
\end{document}